\documentclass{gtart}

\def\ifplaintex{\expandafter\ifx\csname documentclass\endcsname\relax}

\def\gtp{{\mathsurround=0pt\it $\cal G\mskip-2mu$eometry \&\ 
$\cal T\!\!$opology $\cal P\!$ublications}}  

\def\recd{{\small Received:\qua\receiveddate\ifx\reviseddate\relax
\else\qquad Revised:\qua\reviseddate\fi\par}} 

\def\lognumber#1{\def\thelognumber{#1}}
\def\volumenumber#1{\def\thevolumenumber{#1}}
\def\volumeyear#1{\def\thevolumeyear{#1}}
\def\papernumber#1{\def\thepapernumber{#1}}
\def\pagenumbers#1#2{\def\startpage{#1}\def\finishpage{#2}}
\def\published#1{\def\publishdate{#1}}

\def\received#1{\def\receiveddate{#1}}

\def\accepted#1{\def\accepteddate{#1}}

\def\asciiaddress#1{\def\theasciiaddress{#1}}

\let\\\par\let\thelognumber\relax\let\thevolumenumber\relax
\let\thepapernumber\relax\let\thevolumeyear\relax\let\startpage\relax
\let\finishpage\relax\let\publishdate\relax\let\receiveddate\relax
\let\reviseddate\relax\let\accepteddate\relax\let\theasciititle\relax
\let\theasciiauthors\relax\let\theasciiaddress\relax
\let\theasciiabstract\relax

\let\theasciiemail\relax

\ifplaintex
\font\logobig=cmssbx10 scaled 3836
\font\logomed=cmssbx10 scaled 2557
\else
\font\logobig=cmssbx10 scaled 4200
\font\logomed=cmssbx10 scaled 2800
\fi

\long\def\makeagttitle{   
\count0=\startpage
\agt\hfill      
\hbox to 45truept{\vbox to 0pt{\vglue -13truept{\logomed A\kern -.37em{\logobig 
T}\kern -.38em G}\vss}\hss}
\break
{\small Volume \thevolumenumber\ (\thevolumeyear)
\startpage--\finishpage\nl
Published: \publishdate}

\vglue .25truein

{\parskip=0pt\leftskip 0pt plus
1fil\def\\{\par\smallskip}{\Large\bf\thetitle}\par\medskip} \vglue
0.05truein

{\parskip=0pt\leftskip 0pt plus 1fil\def\\{\par}{\sc\theauthors}
\par\medskip}%
 
\vglue 0.03truein 

{\small\leftskip 25truept\rightskip 25truept{\bf Abstract}\stdspace\theabstract

{\bf AMS Classification}\stdspace\theprimaryclass
\ifx\thesecondaryclass\relax\else; \thesecondaryclass\fi\par
{\bf Keywords}\stdspace \thekeywords\par}\vglue 7truept

}   

\ifplaintex
\font\phead=cmsl9 scaled 950
\font\pnum=cmbx10 scaled 913
\font\pfoot=cmsl9 scaled 950
\headline{\vbox to 0pt{\vskip -4.5mm\line{\small\phead\ifnum
\count0=\startpage ISSN 1472-2739 (on-line) 1472-2747 (printed)
\hfill {\pnum\folio}\else\ifodd\count0\def\\{ }%
\ifx\theshorttitle\relax\thetitle\else\theshorttitle\fi\hfill{\pnum\folio}
\else\def\\{ and }{\pnum\folio}\hfill\ifx\theshortauthors\relax\theauthors
\else\theshortauthors\fi\fi\fi}\vss}}
\footline{\vbox to 0pt{\vglue 0mm\line{\small\pfoot\ifnum\count0=\startpage
\copyright\ \gtp\hfill\else
\agt, Volume \thevolumenumber\ (\thevolumeyear)\hfill\fi}\vss}}
\else
\font=cmsl9 scaled 1050
\font\lnum=cmbx10 
\font=cmsl9 scaled 1050
\makeatletter
\def\@oddhead{{\small\lhead\ifnum\count0=\startpage ISSN 1472-2739 
(on-line) 1472-2747 (printed)\hfill {\lnum\number\count0}\else\ifodd\count0
\def\\{ }\ifx\theshorttitle\relax \thetitle \else\theshorttitle\fi\hfill
{\lnum\number\count0}\else\def\\{ and }{\lnum\number\count0}
\hfill\ifx\theshortauthors\relax 
\theauthors\else\theshortauthors\fi\fi\fi}}\def\@evenhead{\@oddhead}
\def\@oddfoot{\small\lfoot\ifnum\count0=\startpage\copyright\ \gtp\hfill\else
\agt, Volume \thevolumenumber\ (\thevolumeyear)\hfill\fi}
\def\@evenfoot{\@oddfoot}
\makeatother
\fi
\let\maketitlepage\makeagttitle
\let\makeshorttitle\maketitlepage
\let\maketitle\maketitlepage

\newwrite\gtoutfile
\long\gdef\makeheadfile{  
{\def\\{, }\def\s{ }
\immediate\openout\gtoutfile head.xxx
\immediate\write\gtoutfile{To: math@arxiv.org}
\immediate\write\gtoutfile{Subject: put OR rep NNNNN:ppppp}
\immediate\write\gtoutfile{--text follows this line--}
\immediate\write\gtoutfile{Proxy-for: \ifx\theasciiauthors\relax
\theauthors\else\theasciiauthors\fi\s<\ifx\theasciiemail\relax\theemail\else\theasciiemail\fi>}
\immediate\write\gtoutfile{\noexpand\\}
\immediate\write\gtoutfile{Authors: \ifx\theasciiauthors\relax
\theauthors\else\theasciiauthors\fi}
{\def\\{ }\immediate\write\gtoutfile{Title: \ifx\theasciititle\relax
\thetitle\else\theasciititle\fi}}
\immediate\write\gtoutfile{Subj-class: GT or SG, GR etc}
\immediate\write\gtoutfile{MSC-class: \theprimaryclass\ifx\thesecondaryclass\relax\else, \thesecondaryclass\fi}
\immediate\write\gtoutfile{Journal-ref: Algebr. Geom. Topol. \thevolumenumber\s
(\thevolumeyear) \startpage-\finishpage}
\immediate\write\gtoutfile{Comments: Published by Algebraic and
Geometric Topology at}
\immediate\write\gtoutfile{\s\s\s  http://www.maths.warwick.ac.uk/agt/AGTVol\thevolumenumber/agt-\thevolumenumber-\thepapernumber.abs.html}
\immediate\write\gtoutfile{\noexpand\\}
\immediate\write\gtoutfile{}
\ifx\theasciiabstract\relax
\immediate\write\gtoutfile{\theabstract}\else
\immediate\write\gtoutfile{\theasciiabstract}\fi
\immediate\write\gtoutfile{}
\immediate\write\gtoutfile{\noexpand\\}
\immediate\write\gtoutfile{}
\immediate\closeout\gtoutfile}}  

\def\maketitlepage{\makeagttitle\makeheadfile}
\let\makeshorttitle\maketitlepage
\let\maketitle\maketitlepage

\def\ifplaintex{\expandafter\ifx\csname documentclass\endcsname\relax}

\def\gtp{{\mathsurround=0pt\it $\cal G\mskip-2mu$eometry \&\ 
$\cal T\!\!$opology $\cal P\!$ublications}}  

\def\recd{{\small Received:\qua\receiveddate\ifx\reviseddate\relax
\else\qquad Revised:\qua\reviseddate\fi\par}} 

\def\lognumber#1{\def\thelognumber{#1}}
\def\volumenumber#1{\def\thevolumenumber{#1}}
\def\volumeyear#1{\def\thevolumeyear{#1}}
\def\papernumber#1{\def\thepapernumber{#1}}
\def\pagenumbers#1#2{\def\startpage{#1}\def\finishpage{#2}}
\def\published#1{\def\publishdate{#1}}

\def\received#1{\def\receiveddate{#1}}

\def\accepted#1{\def\accepteddate{#1}}

\def\asciiaddress#1{\def\theasciiaddress{#1}}

\let\\\par\let\thelognumber\relax\let\thevolumenumber\relax
\let\thepapernumber\relax\let\thevolumeyear\relax\let\startpage\relax
\let\finishpage\relax\let\publishdate\relax\let\receiveddate\relax
\let\reviseddate\relax\let\accepteddate\relax\let\theasciititle\relax
\let\theasciiauthors\relax\let\theasciiaddress\relax
\let\theasciiabstract\relax

\let\theasciiemail\relax

\ifplaintex
\font\logobig=cmssbx10 scaled 3836
\font\logomed=cmssbx10 scaled 2557
\else
\font\logobig=cmssbx10 scaled 4200
\font\logomed=cmssbx10 scaled 2800
\fi

\long\def\makeagttitle{   
\count0=\startpage
\agt\hfill      
\hbox to 45truept{\vbox to 0pt{\vglue -13truept{\logomed A\kern -.37em{\logobig 
T}\kern -.38em G}\vss}\hss}
\break
{\small Volume \thevolumenumber\ (\thevolumeyear)
\startpage--\finishpage\nl
Published: \publishdate}

\vglue .25truein

{\parskip=0pt\leftskip 0pt plus
1fil\def\\{\par\smallskip}{\Large\bf\thetitle}\par\medskip} \vglue
0.05truein

{\parskip=0pt\leftskip 0pt plus 1fil\def\\{\par}{\sc\theauthors}
\par\medskip}%
 
\vglue 0.03truein 

{\small\leftskip 25truept\rightskip 25truept{\bf Abstract}\stdspace\theabstract

{\bf AMS Classification}\stdspace\theprimaryclass
\ifx\thesecondaryclass\relax\else; \thesecondaryclass\fi\par
{\bf Keywords}\stdspace \thekeywords\par}\vglue 7truept

}   

\ifplaintex
\font\phead=cmsl9 scaled 950
\font\pnum=cmbx10 scaled 913
\font\pfoot=cmsl9 scaled 950
\headline{\vbox to 0pt{\vskip -4.5mm\line{\small\phead\ifnum
\count0=\startpage ISSN 1472-2739 (on-line) 1472-2747 (printed)
\hfill {\pnum\folio}\else\ifodd\count0\def\\{ }%
\ifx\theshorttitle\relax\thetitle\else\theshorttitle\fi\hfill{\pnum\folio}
\else\def\\{ and }{\pnum\folio}\hfill\ifx\theshortauthors\relax\theauthors
\else\theshortauthors\fi\fi\fi}\vss}}
\footline{\vbox to 0pt{\vglue 0mm\line{\small\pfoot\ifnum\count0=\startpage
\copyright\ \gtp\hfill\else
\agt, Volume \thevolumenumber\ (\thevolumeyear)\hfill\fi}\vss}}
\else
\font=cmsl9 scaled 1050
\font\lnum=cmbx10 
\font=cmsl9 scaled 1050
\makeatletter
\def\@oddhead{{\small\lhead\ifnum\count0=\startpage ISSN 1472-2739 
(on-line) 1472-2747 (printed)\hfill {\lnum\number\count0}\else\ifodd\count0
\def\\{ }\ifx\theshorttitle\relax \thetitle \else\theshorttitle\fi\hfill
{\lnum\number\count0}\else\def\\{ and }{\lnum\number\count0}
\hfill\ifx\theshortauthors\relax 
\theauthors\else\theshortauthors\fi\fi\fi}}\def\@evenhead{\@oddhead}
\def\@oddfoot{\small\lfoot\ifnum\count0=\startpage\copyright\ \gtp\hfill\else
\agt, Volume \thevolumenumber\ (\thevolumeyear)\hfill\fi}
\def\@evenfoot{\@oddfoot}
\makeatother
\fi
\let\maketitlepage\makeagttitle
\let\makeshorttitle\maketitlepage
\let\maketitle\maketitlepage

\newwrite\gtoutfile
\long\gdef\makeheadfile{  
{\def\\{, }\def\s{ }
\immediate\openout\gtoutfile head.xxx
\immediate\write\gtoutfile{To: math@arxiv.org}
\immediate\write\gtoutfile{Subject: put OR rep NNNNN:ppppp}
\immediate\write\gtoutfile{--text follows this line--}
\immediate\write\gtoutfile{Proxy-for: \ifx\theasciiauthors\relax
\theauthors\else\theasciiauthors\fi\s<\ifx\theasciiemail\relax\theemail\else\theasciiemail\fi>}
\immediate\write\gtoutfile{\noexpand\\}
\immediate\write\gtoutfile{Authors: \ifx\theasciiauthors\relax
\theauthors\else\theasciiauthors\fi}
{\def\\{ }\immediate\write\gtoutfile{Title: \ifx\theasciititle\relax
\thetitle\else\theasciititle\fi}}
\immediate\write\gtoutfile{Subj-class: GT or SG, GR etc}
\immediate\write\gtoutfile{MSC-class: \theprimaryclass\ifx\thesecondaryclass\relax\else, \thesecondaryclass\fi}
\immediate\write\gtoutfile{Journal-ref: Algebr. Geom. Topol. \thevolumenumber\s
(\thevolumeyear) \startpage-\finishpage}
\immediate\write\gtoutfile{Comments: Published by Algebraic and
Geometric Topology at}
\immediate\write\gtoutfile{\s\s\s  http://www.maths.warwick.ac.uk/agt/AGTVol\thevolumenumber/agt-\thevolumenumber-\thepapernumber.abs.html}
\immediate\write\gtoutfile{\noexpand\\}
\immediate\write\gtoutfile{}
\ifx\theasciiabstract\relax
\immediate\write\gtoutfile{\theabstract}\else
\immediate\write\gtoutfile{\theasciiabstract}\fi
\immediate\write\gtoutfile{}
\immediate\write\gtoutfile{\noexpand\\}
\immediate\write\gtoutfile{}
\immediate\closeout\gtoutfile}}  

\def\maketitlepage{\makeagttitle\makeheadfile}
\let\makeshorttitle\maketitlepage
\let\maketitle\maketitlepage

\lognumber{17}
\volumenumber{1}
\volumeyear{2001}
\papernumber{17}
\pagenumbers{349}{368}
\received{8 January 2001}
\accepted{17 May 2001}
\published{2 June 2001}

\let\cal\mathcal
\usepackage{amsmath}
\usepackage{amssymb}
\usepackage{amsthm}
\usepackage{amsgen}
\usepackage{eucal}
\usepackage{graphicx}

\newtheorem{thm}{Theorem}[section]    
\newtheorem{lem}[thm]{Lemma}         
\newtheorem{cor}[thm]{Corolary}
\newtheorem{pro}[thm]{Proposition}

\theoremstyle{definition}
   
\newtheorem*{rem}{Remark}             

\newtheorem*{exa}{Example}
\newtheorem*{claim}{Claim}
\newtheorem*{que}{Question}
\newtheorem*{ques}{Questions}

\begin{document}

\title{A characterization of shortest geodesics on surfaces}
\author{Max Neumann-Coto}

\address{Instituto de Matem\'{a}ticas UNAM\\
Ciudad Universitaria,
M\'{e}xico D.F. 04510, M\'{e}xico}
\asciiaddress{Instituto de Matematicas UNAM,
Ciudad Universitaria,
Mexico D.F. 04510, Mexico}
\email{max@math.unam.mx}

\begin{abstract}
Any finite configuration of curves with minimal intersections on a
surface is a configuration of shortest geodesics for some Riemannian metric
on the surface. The metric can be chosen to make the lengths of these
geodesics equal to the number of intersections along them.
\end{abstract}


\maketitle

If $S$ is a closed surface, with some Riemannian metric, then each
essential curve immersed in $S$ is freely homotopic to a smooth geodesic in 
$S$ which is shortest among all the curves in that homotopy class. So while
closed geodesics represent the critical points of the length in each free
homotopy class, these s{\sl hortest geodesics} represent the absolute minima in
each class. Shortest geodesics have topological properties which are not shared
by all geodesics: Freedman, Hass, Rubinstein and Scott showed in \cite{FHS}
and \cite{HR} that shortest geodesics intersect minimally, i.e., they
have the minimum number of intersections and self-intersections allowed by
their free homotopy classes, unless they factor through coverings of other
shortest geodesics. The curve in figure 1a, for example, represents a
geodesic for some metric on $S$, but it can't represent a shortest geodesic.
On the other hand, when a homotopy class allows different configurations
with minimal intersections, as in figures 1b, 1c and 1d, it seems natural to
ask which ones correspond to shortest geodesics for some Riemannian metric on
$S$.  This question was first considered by Shepard \cite{Sh}.

\begin{figure}[ht!]
\centerline{\includegraphics{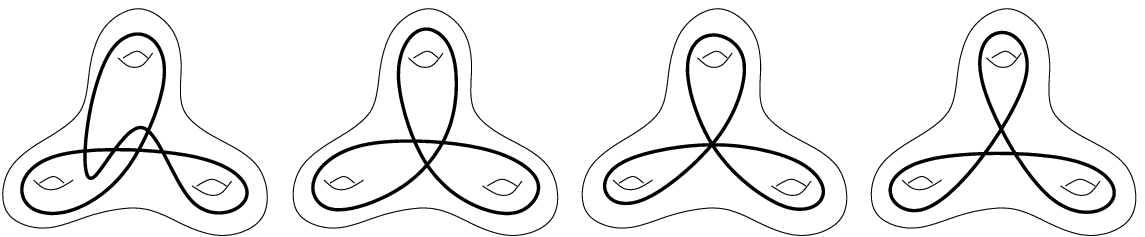}}
\vspace{-2mm}
\centerline{\small a\hspace{1.05in} b\hspace{1.05in} c\hspace{1.05in} d}
\nocolon\caption{}
\end{figure}

In this paper we prove that the shortest geodesics in a surface are
characterized by the minimal intersection property, showing that any finite
configuration of curves with minimal intersections in $S$ is a configuration
of shortest geodesics for some Riemannian metric $g$ on $S$, and also that $%
g $ can be chosen to make the lengths of these geodesics equal to the number
of intersections along them. The proof starts by `blowing up' the metric
outside a regular neighborhood of the curves (an idea introduced by
Bonahon in the context of least area surfaces in 3-manifolds) to transform
the problem into a combinatorial one.
 
The main result implies that all minimal configurations in $S$ can be extended
to contain curves in any other homotopy classes, and gives conditions for
the existence of `absolute' inequalities relating the minimal lengths of
curves in different homotopy classes (inequalities that hold for all
Riemannian metrics on $S$).
In the second part the idea of transmitting cut and paste instructions along a 
homotopy is combined with a result of Hass and Scott \cite{HS} to give a new
proof of the minimal intersection property of \cite{FHS} and to find some absolute 
inequalities involving minimal configurations.

\section{Minimal configurations.}

Two collections of immersed curves in $S$ have
{\sl the same configuration} if there is an ambient isotopy that moves the
image  of one to the other. The curves in a configuration {\sl intersect
minimally} or have {\sl minimal intersections} if they minimize the number of
intersections and self-intersections among all transverse and self-transverse
curves in their free homotopy classes. Following \cite{FHS} and \cite{HR},
the intersections and self-intersections are counted `in the source',
by counting how many curves one crosses when following a curve
all the way around (so multiple intersections are counted with multiplicity and
all the curves in figure 1 have 6 self-intersections). 

As with all geodesics in a surface,
shortest geodesics are transverse and self-transverse, unless they factor
through coverings of other geodesics. Shortest geodesics may not be unique,
and they may not be in general position as they may have points of multiple
intersections.

The results of \cite{FHS} and \cite{HR} can be summarized as follows:

(a)\qua Shortest geodesics intersect minimally, unless they are coverings of
other shortest geodesics.

(b)\qua  If $\alpha$ is an orientation-preserving curve in $S$, the shortest
geodesics representing powers of $\alpha$ always cover a shortest geodesic
representing $\alpha $, but if $\alpha$ is orientation-reversing and there are
2 different shortest geodesics representing $\alpha $, then the shortest
geodesics representing $\alpha^2$ and the odd powers of $\alpha$ do not cover
other shortest geodesics.

According to these results, the image of a collection of shortest geodesics in $S$ is
a configuration of essential curves that intersect transversely and
minimally and do not represent proper powers of any orientation-preserving
class. A finite configuration of essential curves in $S$ with these
properties will be called a {\sl minimal configuration }in $S$.

\thm{Any minimal configuration of curves in
a closed surface $S$ is a configuration of shortest geodesics for some
Riemannian metric g on $S$. g can be chosen so that the length of each curve
in the configuration is equal to the total number of intersections along it.}

\rem{The curves need not be in general position, some
may be homotopic or represent proper powers of an orientation-reversing
class. The curves with no intersections will have length 1.}

\lem{If $c_{1,}c_{2},...,c_{n}$ is a collection of curves immersed
transversely in $S$, and $N$ is a regular neighborhood of $\cup c_{i}$, then
there is a Riemannian metric $g$ on $S$ such that:
\hfill \break
{\rm(a)}\qua Each $c_{i}$ is a geodesic in $S$ and a shortest geodesic in $N$.
\hfill \break
{\rm(b)}\qua All essential curves in $S$ that don't lie entirely in $N$ are longer than
every $c_{i}$. 
\hfill \break
{\rm(c)}\qua The lengths of the arcs of the configuration (the components of 
$\bigcup c_{i}-intersec\-tions$) can be chosen to be any positive numbers.}

\proof The idea is to make the surface look like a
landscape with the curves lying in the bottom of deep and narrow canyons and
surrounded by large mountains.
Start with any Riemannian metric $g_{S}$ on $S$ and regular neighborhoods $%
N^{-}\subset N\subset N^{+}$ of $\bigcup c_{i}$. Since there is a positive
lower bound for the lengths of all essential curves in $S-N^{-}$ and all
arcs running from $S-N$ to $N^{-}$, then by multiplying $g_{S}$ by a
constant $k$ we can make that lower bound larger than the desired lengths of
the curves (this makes the canyons deep and the mountains large).
Now we want to modify the metric inside $N^{+}$. $\bigcup c_{i}$ is a union
of arcs that meet at the multiple points, so $N^{+}$ is the union of
(topological) rectangles and polygons (around the multiple points) as in 
figure 2a. Put a flat metric $g_{N}$ on $N^{+}$ to make these
rectangles and polygons Euclidean, so each $c_{i}$ is a geodesic and no
homotopy of $c_{i}$ within $N^{+}$ reduces its length. Since the lengths of
the rectangles can be chosen independently of each other, we can choose the
length of each arc of the configuration to be any positive number, and the
diameters of the polygons can be taken to be smaller than any prescribed
number $d$ (this makes the canyons long and narrow).

\begin{figure}[ht!]
\centerline{\includegraphics{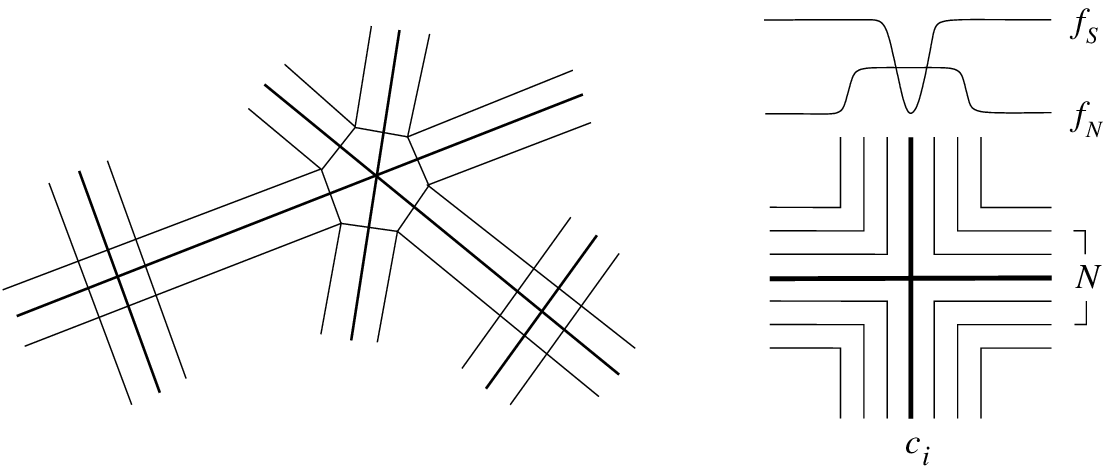}}
\vspace{-4mm}
\centerline{\small a\hspace{2in} b}
\nocolon\caption{}
\end{figure}

Let $g=f_{S}g_{S}+f_{N}g_{N}$ (in terms of the first fundamental forms),
where $f_{S}$ and $f_{N}$ are smooth scalar functions such that:

$f_{N}=1$ on $N$, $f_{N}>0$ on $N^{+}-N$ and $f_{N}=0$ on $S-N^{+}$.

$f_{S}=0$ on $\bigcup c_{i}$ , $f_{S}>0$ on $N^{-}-\bigcup c_{i}$
and $f_{S}=k$ on $S-N^{-}.$

See figure 2b. As $g\geq kg_{S}$ on $S-N^{-}$, with the metric $g$
any essential curve in $S-N^{-}$ and any arc that crosses from $N^{-}$ to $%
S-N$ is longer than $c_{i}$, so any essential curve which doesn't lie in $N$
is longer than $c_{i}$. As $c_{i}$ was a shortest geodesics in $N$ with the
metric $g_{N}$ , and $g=g_{N}$ on $\bigcup c_{i}$ but $g>g_{N}$ on $%
N-\bigcup c_{i}$, the metric $g$ makes every $c_{i}$ a geodesic in $S$ and
shortest geodesic in $N$. \endproof

Each essential curve $c$ immersed in $N$ is freely homotopic in $N$ to a
{\sl polygonal} curve $p$ made of arcs of the configuration (maybe repeated)
and we may assume that $p$ is reduced in the sense that no arc is followed
immediately by the same arc in the opposite direction.

\lem{The metric $g$ can be defined so that if
a (reduced) polygonal curve $p$ is longer than $c_{i}$, then all the curves
homotopic to $p$ in $N$ are also longer than $c_{i}$.}

\proof This is accomplished by choosing $d$ small
(narrow canyons). If $p$ has corners (i.e., if $p$ is not one of the
$c_{i}$'s) then its length can be reduced by rounding the corners, but no
homotopy within $N$ can reduce its length by more than $d$ multiplied by the
number of corners of $p$ (this is clear for the metric $g_{N}$, and
$g=g_{N}$ along  $p$ and $g\geq g_{N}$ elsewhere in $N$), and the number of
corners in $p$ is bounded above by a linear function of its length. 
So if $c$ is homotopic to $p$ in $N$ then $\frac{length(c)}{length(p)}>1-ld$,
where $l$ is the coefficient of the linear function, so by taking $d$ small
enough we can make this ratio as close to 1 as we want. But the set of lengths
of the polygonal curves in $N$ is discrete (because it is contained in the set
of positive linear combinations of the lengths of the arcs), so $\frac{%
length(c_{i})}{length(p)}<m<1$ for all $c_{i}$ 's and all longer $p$'s. So
by taking $d$ small we can make $\frac{length(c)}{length(p)}>\frac{%
length(c_{i})}{length(p)}$ for every $p$ longer than $c_{i}$.\endproof

The previous lemmas have no minimal intersection hypothesis: the metric $g$
makes each $c_{i}$ a geodesic in $S$, but not necessarily a shortest
geodesic, because nonhomotopic curves in $N$ may be homotopic in $S$. The $%
c_{i}$'s are shortest geodesics for some metric $g$ on $S$ if and only if
the lengths of the arcs of the configuration can be chosen so that all
homotopic $c_{i}$'s have the same length and all polygonal curves homotopic
to $c_{i}$ are longer than $c_{i}$.

Now let $c_{1,}c_{2},...,c_{n}$ be a collection of curves with minimal
intersection and self-intersection in $S$. In order to choose the lengths of
the arcs of the configuration, take a collection of {\sl measuring curves}
$\mu _{1,}\mu _{2},...,\mu _{m}$ in general position with respect to the $%
c_{i}$'s, assign to each $\mu _{j}$ a positive {\sl width} $w_{j}$, and
define the length of each arc of the configuration as the sum of the widths
of the curves $\mu _{j}$ that meet the arc. As the arcs of the configuration
must have positive length, we need a collection $\left\{ \mu _{j}\right\} $
whose union meets all the arcs.

Let's say that a measuring collection for $\left\{ c_{i}\right\} $ is 
{\sl good} if it intersects each $c_{i}$ minimally but does not intersect
any polygonal curve homotopic to some $c_{i}$ minimally.

\lem{If $\left\{ \mu _{j}\right\} $ is a good
measuring collection, then for any choice of widths the assigned lengths
make the $c_{i}$'s shortest geodesics for a Riemannian metric in $S$}.

\proof As the $\mu _{j}$'s intersect $c_{i}$ minimally,
if a polygonal curve $p$ is homotopic to $c_{i}$, then each $\mu _{j}$ must
intersect $p$ at least as many times as it intersects $c_{i}$, so $p$ is at
least as long as $c_{i}$, and it is longer than $c_{i}$ if and only if the
total number of intersections of the $\mu _{j}$'s with $p$ is larger, i.e.,
if some $\mu _{j}$ does not intersect $p$ minimally. This is clearly
independent of the choice of widths. Now apply lemmas 1.2 and 1.3. \endproof

\rem{Notice that if a good measuring collection is extended
in any way (by adding curves that intersect the $c_{i}$'s minimally) then
the resulting measuring collection is good.}

\bigskip

{\bf Construction of good measuring collections}

\medskip

Let $\left\{ c_{i}\right\} $ be a configuration of essential curves in a
surface $S$. 
If $\chi (S)\leq 0$, the universal covering of $S$ is a
plane $\widetilde{S}$, and the cyclic coverings $S^{\alpha }$ of $S$ corresponding
to the subgroups generated by elements $\alpha $ of $\pi _{1}(S)$ are annuli
or Moebius bands (depending on whether $\alpha$ is orientation preserving
or orientation-reversing). 
So the preimage of $\left\{ c_{i}\right\}$ in $\widetilde{S}$ is an infinite 
configuration of topological lines, 
while the preimage of $\left\{c_{i}\right\}$ in $S^{\alpha}$ is a configuration 
of lines and curves (the liftings of the $c_{i}$'s representing
powers of $\alpha$, if any). The curves in $S^{\alpha}$ will be denoted by
$c_{i}^{\alpha }$ and the lines in $S^{\alpha}$ or $\tilde{S}$ by
$\widetilde{c}_{i}$.

According to \cite{FHS} and \cite{HR}, $\left\{ c_{i}\right\}$ is a minimal
configuration in $S$ if and only if for each $S^{\alpha}$ 
the curves $c_{i}^{\alpha }$ intersect minimally and intersect the lines
$\widetilde{c}_{i}$ minimally, that is:

(a)\qua If $\alpha $ is orientation-preserving then the curves representing $\alpha$
in $S^{\alpha }$ are embedded and disjoint, and
intersect each line in at most 1 point.

(b)\qua  If $\alpha $ is orientation-reversing then the curves representing $\alpha^{r}$, 
$r$ odd, have $r-1$ self-intersections and intersect the curves
representing $\alpha ^{s}$ (s odd, $s\geqslant r$) in $r$ points.
The curves representing $\alpha ^{2}$ are embedded and
disjoint from all the other curves representing powers of $\alpha$.
A curve representing $\alpha ^{r}$ ($r=2$ or odd) intersects a line 
that crosses $S^{\alpha }$ in $r$ points.

\medskip

{\bf Case 1}\qua {\sl All $c_{i}$'s are primitive and
orientation-preserving.}

A natural candidate for a good measuring collection consists of a pair of 
'parallel' curves $\mu _{i+}$ and $\mu _{i-}$ for each $c_{i}$, one to
the right and one to the left of $c_{i}$ and sufficiently close so that the
immersed annulus determined by $\mu _{i+}$ and $\mu _{i-}$ intersects the
curves of the configuration along arcs that cross the annulus, and the only
multiple points of the configuration inside the annulus are the ones along $%
c_{i}$. So $c_{i}$ and $\mu _{i+}$ intersect each $c_{j}$ the same number of
times (the arcs of intersection between the curves and the annulus give a
one to one correspondence between the intersections along $c_{i}$ and the
intersections along $\mu _{i+}$) and so $\mu _{i+}$ intersects each $c_{j}$
minimally. By construction $\left\{ \mu _{i+},\mu _{i-}\right\} $ meets all
the arcs of the configuration. Notice that taking all the widths equal to $%
\frac{1}{2}$ makes the length of each $c_{i}$ equal to the number of
intersections of the configuration along $c_{i}$ (counted with
multiplicity).

Define the {\sl distance} between two lines in the configuration $\left\{\widetilde{c}_{i}\right\}$
in $\widetilde{S}$ as the minimum number of complementary regions 
that one has to cross to go from one  line to the other (so
the distance is 0 iff the lines meet). We will say that 2 -not necessarily
different- curves $c_{i}$ and $c_{j}$ in $S$ are {\sl close neighbors} if two
of their preimages $\widetilde{c}_{i}$ and $\widetilde{c}_{j}$ are at
distance 1 in $\widetilde{S}$.

\claim In case 1, {$\left\{ \mu _{j+},\mu _{j-}\right\}$ is
a good measuring collection for $\left\{{c_{i}}\right\}$
if and only if every $c_{i}$ has close neighbors on both sides.}

\proof Observe that a polygonal curve $p$ homotopic to $%
c_{i}$ intersects $\left\{ \mu _{i+},\mu _{i-}\right\} $ minimally if and
only if in the corresponding covering $S^{\alpha }$, the curves $p^{\alpha }$
and $c_{i}^{\alpha }$ intersect the same measuring curves and lines ($\mu
_{i+}^{\alpha }$, $\mu _{i-}^{\alpha }$, $\widetilde{\mu }_{j+}$, $%
\widetilde{\mu }_{j-}$) and do so the same number of times (once in this
case).
So $p^{\alpha }$ cannot cross or touch $c_{i}^{\alpha }$ or any other curve $%
c_{j}^{\alpha }$ (because $p^{\alpha }$ would intersect $\mu _{i+}^{\alpha } 
$, $\mu _{i-}^{\alpha }$, $\mu _{j+}^{\alpha }$ or $\mu _{j-}^{\alpha }$),
and the annulus bounded by $c_{i}^{\alpha }$ and $p^{\alpha }$ must
intersect the lines $\widetilde{c}_{j}$ along arcs that cross the annulus
(if a line $\widetilde{c}_{j}$ touches this annulus at one point or
intersects it along an arc that starts and ends in $p^{\alpha }$, then $%
p^{\alpha }$ intersects one of the lines $\widetilde{\mu }_{j+}$ or $%
\widetilde{\mu }_{j-}$ twice). In particular $p^{\alpha }$ must be made
exclusively of arcs of lines $\widetilde{c}_{j}$ that cross $c_{i}^{\alpha }$.

If all the $c_{i}$'s representing $\alpha $ have close neighbors on both
sides, then the curves $c_{i}^{\alpha }$ in $S^{\alpha }$ are just one
complementary region away from other $c_{j}^{\alpha }$'s or from lines $%
\widetilde{c}_{j}$ that don't meet $c_{i}^{\alpha }$. So any polygonal curve 
$p^{\alpha }$ must cross or at least touch one of these curves or lines, and
so $p$ cannot have minimal intersections with $\left\{ \mu _{j+},\mu
_{j-}\right\} $.

\begin{figure}[ht!]
\centerline{\includegraphics{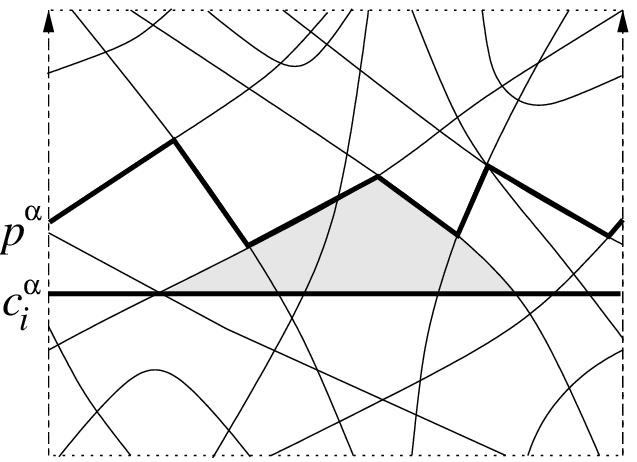}}
\nocolon\caption{}
\end{figure}

Now suppose that some $c_{i}$ representing $\alpha $ doesn't have close
neighbors on one side. Then for each complementary region $R$ on one side of 
$c_{i}^{\alpha }$, the lines adjacent to $R$ must intersect $c_{i}^{\alpha }$%
. These lines determine triangles with base in $c_{i}^{\alpha }$ that
contain $R$, and the widest of these triangles (the one with maximal base in 
$c_{i}^{\alpha }$) is crossed by the lines $\widetilde{c}_{j}$ along arcs
that meet the base of the triangle. See figure 3. The union of these wide
triangles on one side of $c_{i}^{\alpha }$ is an annulus whose boundaries
are $\widetilde{c}_{i}$ and a polygonal curve $p^{\alpha }$, and the lines $%
\widetilde{c}_{j}$ can intersect this annulus only along arcs that cross the
annulus, so $p^{\alpha }$ projects to a polygonal curve $p$ homotopic to $%
c_{i}$ that intersects $\left\{ \mu _{i+},\mu _{i-}\right\} $ minimally, and 
$p$ is the nearest polygonal curve with this property.\endproof

In the configuration of figure 4a the curve $c_{1}$ has close
neighbors on both sides but the curve $c_{2}$ doesn't. So a metric $g$ that
makes each arc of the configuration of length 1 can make $c_{1}$ (but not $%
c_{2}$) a shortest geodesic in the surface.

\begin{figure}[ht!]
\centerline{\includegraphics{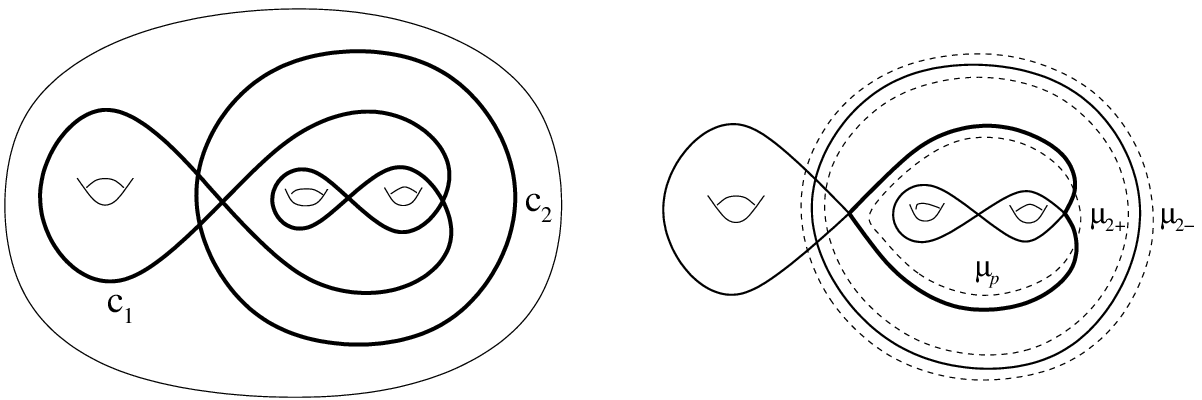}}
\vspace{-2mm}
\centerline{\small a\hspace{2in} b}
\nocolon\caption{}
\end{figure}

Now we want to extend $\left\{ \mu _{i+},\mu _{i-}\right\} $ to a good
measuring collection in the case that some $c_{i}$'s don't have close
neighbors. The idea is given in figure 4b: 
if there is a polygonal curve $p$ homotopic to $c_{i}$ that
intersects $\left\{ \mu _{i+},\mu _{i-}\right\} $ minimally, take a
measuring curve $\mu _{p}$ that runs ``quasiparallel'' to $p$ crossing each
edge of $p$ once, so its lifting to $S^{\alpha }$ looks like in figure 5,
making it sufficiently close so that each $c_{j}$ intersects the singular
annulus determined by $c_{i}$ and $\mu _{p}$ along arcs that cross the
annulus.

\begin{figure}[ht!]
\centerline{\includegraphics[width=2in]{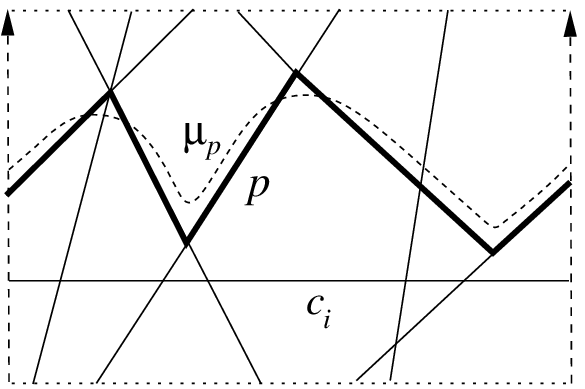}}
\nocolon\caption{}
\end{figure}

To see that $\mu _{p}$ intersects each $c_{j}$ minimally, it is enough to
show that $\mu _{p}$ and $c_{i}$ intersect each $c_{j}$ the same number of
times. This happens because in $S^{\alpha }$ each line $\widetilde{c}_{j}$
intersects the annulus determined by $c_{i}^{\alpha }$ and $\mu _{p}^{\alpha
}$ along arcs that cross the annulus (an arc of intersection of $\widetilde{c%
}_{j}$ with the annulus cannot start and end in $\mu _{p}^{\alpha }$,
because then $\widetilde{c}_{j}$ would cross $p^{\alpha }$ twice or it would
touch it at one point). So by adding a measuring curve $\mu _{p}$ for each
short polygonal curve $p$ we can extend $\left\{ \mu _{i+},\mu _{i-}\right\} 
$ to a good measuring collection. As only finitely many polygonal curves
homotopic to $c_{i}$ can intersect $\left\{ \mu _{i+},\mu _{i-}\right\} $
minimally (because the number of arcs in such polygonals is bounded above by
the number of intersections along $c_{i}$), we are done.

\medskip

{\bf Case 2}\qua {\sl All orientation-reversing $c_{i}$'s are
primitive and no two of them are homotopic}

Choose the measuring curves corresponding to the orientation-preserving
curves as in case 1. The orientation-reversing $c_{i}$'s are one sided, so
instead of two parallel curves $\mu _{i+}$ and $\mu _{i-}$ there is a single
curve $\mu _{i\pm }$ homotopic to $c_{i}^{2}$ that runs on ``both sides''
of $c_{i}$.
To see that $\mu _{i\pm }$ intersects each $c_{j}$ minimally, look at the
covering $S^{\alpha }$ of $S$ corresponding to the class $\alpha $
represented by $c_{i}$. The only closed curve in $S^{\alpha }$ is $%
c_{i}^{\alpha }$, which by construction does not meet $\mu _{i\pm }^{\alpha
} $, and the lines $\widetilde{c}_{j}$ that cross $S^{\alpha }$ intersect $%
c_{i}^{\alpha }$ at a single point, so they must intersect $\mu _{i\pm
}^{\alpha }$ at exactly two points.
Now if $p$ is any polygonal curve homotopic to $c_{i}$, then (as $c_{i}$ is
one sided) $p^{\alpha }$ must cross $c_{i}^{\alpha }$, so $p^{\alpha }$ must
intersect $\mu _{i\pm }^{\alpha }$, and so $\mu _{i\pm }$ doesn't intersect 
$p$ minimally. Therefore this measuring collection is already good.

Observe that if an orientation-reversing $c_{i}$ is nonprimitive, or is
homotopic to another $c_{j}$, then $\mu _{i\pm }$ doesn't have minimal
intersection with $c_{i}$ (or $c_{j}$), and therefore $\mu _{i\pm }$ cannot
be used as a measuring curve.

\medskip

{\bf Case 3}\qua {\sl $S$ is a projective plane}

All $c_{i}$'s are homotopic to the unique nontrivial element of $\pi _{1}(S)$, 
so they are embedded and intersect each other in 1 point. For each $c_{i}$
take a collection of measuring curves $\mu _{ix}$ each made of an arc that
runs parallel to $c_{i}$ all the way around and a small arc that crosses $%
c_{i}$ at one point, as in figure 6a. 

Make $\mu _{ix}$ sufficiently close to 
$c_{i}$ so that the other $c_{i}$'s intersect the singular strip determined
by $c_{i}$\ and $\mu _{ix}$ along arcs that cross it from $c_{i}$ to $\mu
_{ix}$ (so $\mu _{ix}$ will intersect each $c_{j}$ once) and the only
multiple points of the configuration inside the band are the ones along $%
c_{i}$. Take two $\mu _{ix}$ for each arc of $c_{i}$, one crossing the arc
in each direction.
Now if $p$ is any polygonal curve homotopic to $c_{i}$ and $p$ contains arcs
of $c_{j}$, then some $\mu _{jx}$ crosses a corner of $p$ twice (see figure
6b), so $p$ doesn't have minimal intersection with that $\mu _{jx}$.

\begin{figure}[ht!]
\centerline{\includegraphics{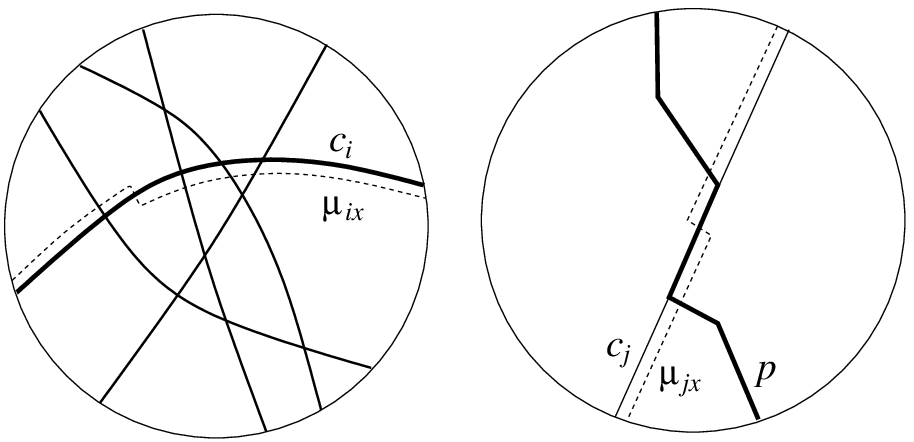}}
\vspace{-2mm}
\centerline{\small a\hspace{2in} b}
\nocolon\caption{}
\end{figure}

{\bf Case 4}\qua {\sl All orientation-reversing $c_{i}$'s are
primitive, but some are homotopic}

Choose the measuring curves for the orientation-preserving $c_{i}$'s as in
case 1, and those for the orientation-reversing $c_{i}$'s that are not
homotopic to other $c_{j}$'s as in case 2.

Now consider an orientation-reversing class $\alpha $ represented by 2 or
more $c_{i}$'s. These $c_{i}$'s lift to curves $c_{i}^{\alpha }$ in the
Moebius band $S^{\alpha }$ that are embedded and intersect each other in 1
point.
For each of these $c_{i}$'s take a collection of measuring curves $\mu _{ix}$
as in case 3, each one made of an arc that runs parallel to $c_{i}$ all the
way around and a small arc that crosses $c_{i}$ at one point, so $\mu _{ix}$
lifts to a curve $\mu _{ix}^{\alpha }$ in $S^{\alpha }$ that intersects $%
c_{i}^{\alpha }$ in exactly one point as in figure 7a. Take again one $\mu
_{ix}$ crossing each arc of $c_{i}$ in each direction. As $\mu _{ix} $ and $%
c_{i}$ intersect each $c_{j}\neq c_{i}$ the same number of times, then $\mu
_{ix}$ intersects each $c_{j}\neq c_{i}$ minimally, and as $\mu _{ix}$
intersects $c_{i}$ one more time than $c_{i}$ intersects itself, then $\mu
_{ix}$ also intersects $c_{i}$ minimally.

\begin{figure}[ht!]
\centerline{\includegraphics{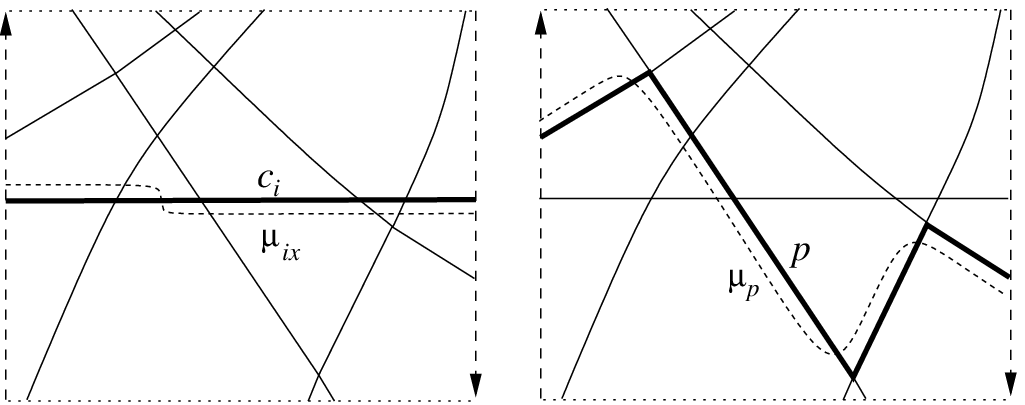}}
\vspace{-1mm}
\centerline{\small a\hspace{2in} b}
\nocolon\caption{}
\end{figure}

Observe that a polygonal curve $p$ that intersects these $\mu _{ix}$'s
minimally must be made exclusively of arcs of orientation-preserving $c_{j}$'s, 
because if $p$ contains an arc of some orientation-reversing $c_{j}$
then one of the $\mu _{jx}$'s crosses a corner of $p$ twice. And one can
show as in case 1 that such $p$ has minimal intersection with the curves $%
\mu _{i+}$, $\mu _{i-}$ and $\mu _{ix}$ if and only if $p^{\alpha }$
intersects $c_{i}^{\alpha }$ and every curve $c_{j}^{\alpha }$ homotopic to $%
c_{i}^{\alpha }$ in exactly one point, and each line $\widetilde{c}_{j}$
that intersects the singular annulus determined by $c_{i}^{\alpha }$ and $%
p^{\alpha }$ does so along one arc that crosses the annulus from $p^{\alpha
} $ to $c_{i}^{\alpha }$.

If there is a short polygonal curve $p$ homotopic to $\alpha $, take a
measuring curve $\mu _{p}$ that runs ``quasiparallel'' to $p$ crossing each
edge of $p$ once, so $\mu _{p}$ lifts to a curve $\mu _{p}^{\alpha }$ in $%
S^{\alpha }$ that looks like in figure 7b. To see that $\mu _{p}$ intersects
each $c_{j}$ minimally, observe that $\mu _{p}^{\alpha }$ intersects each
curve $c_{j}^{\alpha }$ once (otherwise $p^{\alpha }$ would intersect $%
c_{j}^{\alpha }$ more than once) and that each $\widetilde{c}_{j}$
intersects the singular annulus determined by $c_{i}^{\alpha }$ and $\mu
_{p}^{\alpha }$ along arcs that cross it from $c_{i}^{\alpha }$ to $\mu
_{p}^{\alpha }$ (an arc of intersection cannot start and end in $\mu _{p}$,
because then $\widetilde{c}_{j}$ would intersect the singular annulus
determined by $c_{i}^{\alpha }$ and $p^{\alpha }$ in an arc that starts and
ends in $p^{\alpha }$).

By construction the number of intersections between $\mu _{p}^{\alpha }$ and 
$p^{\alpha }$ is equal to the number of corners of $p$, so $\mu _{p}$
intersects $p$ minimally only when $p$ has one corner.
To deal with these short polygonal curves with only one corner, we need an
extra measuring curve $\mu _{\alpha \pm }$ whose lifting to $S^{\alpha }$
runs parallel to the boundary of the region $V\alpha $ determined by all the
curves $c_{i}^{\alpha }$, as in figure 8a, so $\mu _{\alpha \pm }$ is
homotopic to $\alpha ^{2}$. $\mu _{\alpha \pm }$ intersects each $c_{j}$
minimally because the curves $c_{i}^{\alpha }$ are contained in the Moebius
band bounded by $\mu _{\alpha \pm }$, and if a line $\widetilde{c}_{j}$
intersects this Moebius band along a nonessential arc then $\widetilde{c}_{j}
$ intersects some $c_{i}^{\alpha }$ in two points.

\begin{figure}[ht!]
\centerline{\includegraphics{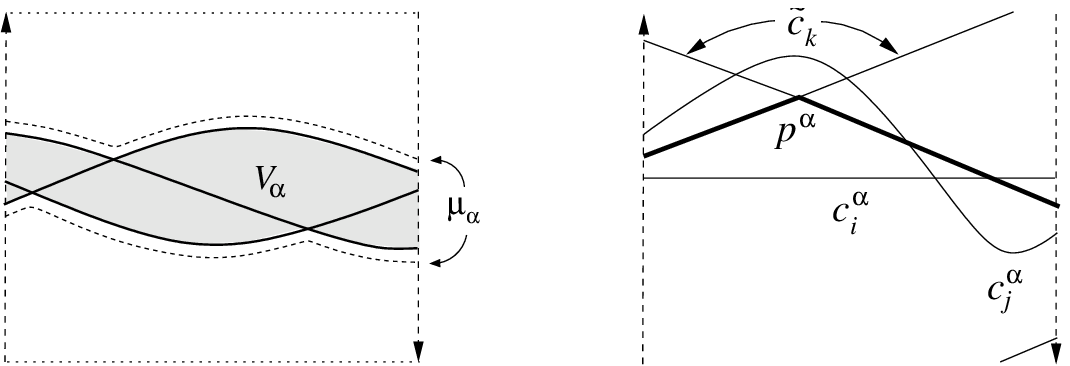}}
\vspace{-2mm}
\centerline{\small a\hspace{2.5in} b}
\nocolon\caption{}
\end{figure}

We claim that if $p$ is a short polygonal curve with one corner then $%
p^{\alpha }$ cannot be contained in $V\alpha $, so $p^{\alpha }$ intersects $%
\mu _{\alpha \pm }$ and so $p$ doesn't intersect $\mu _{\alpha \pm }$
minimally. If $p^{\alpha }$ were contained in $V\alpha $ then its corner
would be in the region determined by two curves $c_{i}^{\alpha }$ and $%
c_{j}^{\alpha }$. As $p^{\alpha }$ is made of an arc of a line $\widetilde{c}%
_{k}$ that starts and ends at the corner, $\widetilde{c}_{k}$ would have to
cross $c_{i}^{\alpha }$ or $c_{j}^{\alpha }$ twice (see figure 8b)
contradicting the fact that $c_{k}$ intersects $c_{i}$ and $c_{j}$ minimally.

\medskip

{\bf Case 5}\qua {\sl Some orientation-reversing $c_{i}$'s are
nonprimitive}

Let $1\leq r_{1}<r_{2}<...<r_{n}$ be the odd powers of a primitive
orientation-reversing class $\alpha $ represented by $c_{i}$'s in the
configuration. Each of these $c_{i}$'s lifts to an immersed curve $%
c_{i}^{\alpha }$ in the Moebius band $S^{\alpha }$.
For each of these $c_{i}$'s take a collection of measuring curves $\mu _{ix}$
as in case 4, each made of an arc that runs parallel to $c_{i}$ all the way
around and a small arc that crosses $c_{i}$ at one point. So each $\mu _{ix}$
intersects every $c_{j}$ minimally and every polygonal curve $p$ that
intersects these $\mu _{ix}$'s minimally is made of arcs of
orientation-preserving $c_{j}$'s
For each polygonal curve $p$ homotopic to $c_{i}$ that intersects these $\mu
_{ix}$'s minimally, take a measuring curve $\mu _{p}$ that runs
quasiparallel to $p$ crossing each edge of $p$ once so, as in case 4, $\mu
_{p}$ intersects every $c_{j}$ minimally, but $\mu _{p}$ intersects $p$
minimally only when $p$ has one corner.

To deal with these short polygonal curves with one corner representing $%
\alpha ^{r_{k}}$, we need to add an extra measuring curve $\mu _{\alpha \pm
} $ homotopic to $\alpha ^{2}$ and measuring curves $\mu _{\alpha k}$
homotopic to $\alpha ^{r_{k+1}}$for each $k<n$.
One can show as in case 4 that the polygonal curves with one corner
representing $\alpha ^{r_{k}}$ cannot be contained in the region $V\alpha
^{r_{k}}$ of $S^{\alpha }$ determined by the images of all the $%
c_{i}^{\alpha }$'s representing $\alpha ^{r_{k}}$. The minimal intersection
of the curves in $S^{\alpha }$ implies that all the curves representing $%
\alpha ^{r_{k}}$ are contained in the region determined by each curve
representing $\alpha ^{r_{k+1}}$, so $V\alpha ^{r_{k}}\subset V\alpha
^{r_{k+1}}$, each curve representing $\alpha ^{r_{k+1}}$ intersects $V\alpha
^{r_{k}}$ along one arc, and each line that intersects $V\alpha ^{r_{k}}$
does so along one essential arc.

Let $\mu _{\alpha \pm }$ be a curve whose lifting to $S^{\alpha }$ runs
parallel to the boundary of $V\alpha ^{r_{n}}$, so $\mu _{\alpha \pm }$ is
homotopic to $\alpha ^{2}$. Then $\mu _{\alpha \pm }$ intersects every $%
c_{j} $ minimally, but any polygonal curve representing an odd power of $%
\alpha $ that intersects $\mu _{\alpha \pm }$ minimally must be contained in 
$V\alpha ^{r_{n}}$.
Now for each $k<n$, choose a curve $c_{i}^{\alpha }$ representing $\alpha
^{r_{k+1}}$ which is closest to $V\alpha ^{r_{k}}$ in the sense that the
region determined by its image does not contain any other $c_{j}^{\alpha }$
representing $\alpha ^{r_{k+1}}$. Let $\mu _{\alpha k}$ be a curve whose
lifting to $S^{\alpha }$ runs parallel to the arc $c_{i}^{\alpha }\cap
V\alpha _{r_{k}}$ and then runs around the boundary of $V\alpha _{r_{k}}$
enough times to complete a curve homotopic to $\alpha ^{r_{k+1}}$.

Figure 9a shows a lifting of $\mu _{\alpha k}$ to $S^{\alpha ^{r_{k+1}}}$.
To prove that $\mu _{\alpha k}$ intersects each $c_{j}$ minimally, 
it is enough to show that in the covering $S^{\alpha ^{r_{k+1}}}$ 
the preimages of $c_{j}$ intersect the region determined by the 
liftings of $c_{i}$ and $\mu _{\alpha k}$ along arcs that
cross that region.
An arc of intersection $a$ that didn't cross that region would look as
in figure 9b, but this arc cannot belong to a line $\widetilde{c} _{j}$
because then $\widetilde{c}_{j}$ would intersect $V\alpha ^{r_{k}}$ in at
least two arcs, and it cannot belong to a curve $c_{j}^{\alpha }$ representing 
$\alpha ^{r_{k}}$ or a smaller power of $\alpha $ because these curves are
contained in $V\alpha ^{r_{k}}$. So the arc $a$ must belong to a curve
$c_{j}^{\alpha }$ representing some larger power of $\alpha $, and so
$c_{j}^{\alpha }=a \cup a'$, where $a'$ is an arc  in $V\alpha ^{r_{k}}$. So
$c_{j}^{\alpha }$ lies in the region determined by $c_{i}^{\alpha }$, but by
the choice of $c_{i}^{\alpha }$ no curve representing $V\alpha ^{r_{k+1}}$ or
a larger power of $\alpha $ can be contained in this region.

Now if $p$ is a polygonal curve with one corner representing $\alpha ^{r_{k}}$ 
then its lifting to $S^{\alpha }$ is not contained in $V\alpha ^{r_{k}}$,
so it is not contained in the region determined by $\mu _{\alpha k}$, so $p$
does not intersect $\mu _{\alpha k}$ minimally.

\begin{figure}[ht!]
\centerline{\includegraphics{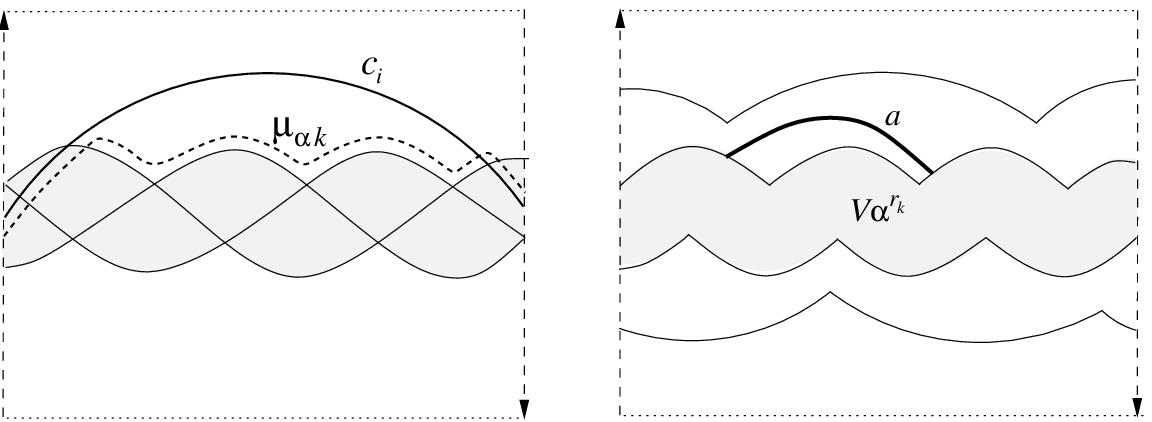}}
\vspace{-2mm}
\centerline{\small a\hspace{2.5in} b}
\nocolon\caption{}
\end{figure}

\bigskip

{\bf Choice of widths}

\medskip

The measuring collection for $\left\{ c_{i}\right\} $ constructed above is
made of curves homotopic to some $c_{i}$ ($\mu _{i+}$, $\mu _{i-}$, $\mu _{i\times}$, 
$\mu _{p}$ and $\mu _{\alpha k}$) or the square of some
primitive orientation-reversing class $\alpha $ ($\mu _{i\pm }$ and $\mu
_{\alpha \pm }$). The choice of widths to prove the second part of the
theorem is not obvious because the minimum number of self-intersections of an
orientation-reversing $c_{i}$ differs from the minimum number of intersections
between $c_{i}$ and a homotopic curve by 1. The condition that the
length of each $c_{i}$ in the configuration must be equal to the number of
intersections along it gives a system of linear equations on the widths of
the measuring curves that has a unique solution for the sums of widths of the
measuring curves in each homotopy class:

(a)\qua Make the sum of the widths of the measuring curves ($\mu _{i+}$, $\mu
_{i-}$ and $\mu _{p}$'s) in each orientation-preserving class\ equal to the
number of $c_{i}$'s in that class.

(b)\qua If $\alpha $ is a primitive orientation-reversing class, and $1$ $\leq
r_{1}<r_{2}<...$ are the odd powers of $\alpha $ represented by some $c_{i}$
's, make the sum of the widths of the measuring curves ($\mu _{ix}$, $\mu
_{p}$ and $\mu _{\alpha k-1}$'s) representing $\alpha ^{r_{k}}$, $k>1$,
equal to the number of $c_{i}$'s representing that class, but for the
measuring curves representing $\alpha ^{r_{1}}$ make the sum of their widths 
$\frac{1}{r_{1}}${\sl \ }units less than the number of $c_{i}$'s
representing that class. Finally, make the width of each measuring curve ($%
\mu _{i\pm }$ and $\mu _{\alpha \pm }$) representing $\alpha ^{2}$\ equal to 
$\frac{1}{2}$.

A problem arises when $r_{1}=1$ and only 1 curve $c_{i}$ represents $\alpha
^{r_{1}}$, because then the\ sum of the widths of the measuring curves
homotopic to $\alpha $ is 0, which means that these measuring curves cannot
be used, and the rest of the measuring collection may not be good. This can
be arranged by replacing the measuring curves representing $\alpha $ by
suitable curves representing $\alpha ^{r_{2}}$ as follows:

\begin{figure}[ht!]
\centerline{\includegraphics{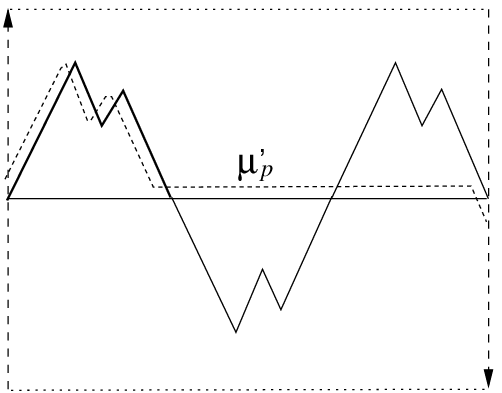}}
\nocolon\caption{}
\end{figure}

Trade each curve $\mu _{ix}$ made of an arc that goes once around $c_{i}$
and a small arc that crosses $c_{i}$ at one point, for a curve $\mu
_{ix}^{\prime }$ made of an arc that goes $r_{2}$ times around $c_{i}$ and
the small arc. And trade each $\mu _{p}$ representing $\alpha $ for a curve $%
\mu _{p}^{\prime }$ obtained by replacing the small arc of $\mu _{p}$ that
crosses $c_{i}$ by an arc that goes $r_{2}-1$ times around $c_{i}$ so it now
represents $\alpha ^{r_{2}}$ (figure 10 shows the lifting of $\mu
_{p}^{\prime }$ to $S^{\alpha ^{r_{2}}}$). It is not hard to see that $\mu
_{ix}^{\prime }$ and $\mu _{p}^{\prime }$ intersect each $c_{j}$ minimally,
but intersect nonminimally all the polygonal curves that intersect $\mu _{ix}
$ or $\mu _{p}$ nonminimally. This proves the second part of the theorem.
\endproof 

\begin{figure}[ht!]
\centerline{\includegraphics{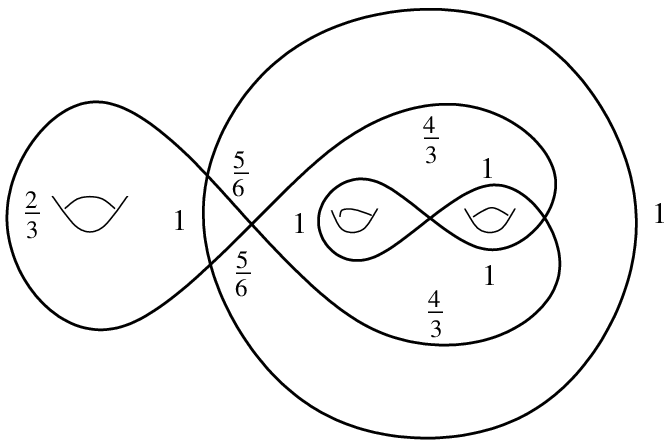}}
\nocolon\caption{}
\end{figure}

Figure 11 shows the lengths of the arcs in the configuration in figure 4 resulting from
making $\mu _{i+}$ and $\mu _{1-}$ of width $\frac{1}{2}$ and $\mu _{2-}$, $%
\mu _{2+}$ and $\mu _{p}$ of width $\frac{1}{3}$.

\rem{Theorem 1.1 clearly holds for nonclosed
surfaces, provided that the curves don't meet the boundary.
It also works for minimal configurations of properly immersed
curves and arcs in a surface with boundary, considering either minimal
configurations with the endpoints of the arcs fixed or free to move along $%
\partial S$ (one just needs to use measuring arcs analogous to the measuring
curves).}

Theorem 1.1 contrasts with the examples of Hass and Scott \cite{HS2} of
minimal configurations of primitive and nonhomotopic curves in a surface
which are not configurations of geodesics for any metric of negative curvature
on the surface. These configurations, however, can be realized by metrics of 
non-positive curvature.

\ques{Which configurations of primitive curves in a surface are configurations
of shortest geodesics for metrics of negative curvature? and for metrics of non
positive curvature?}

The second part of theorem 1.1 is only significant for configurations
containing more than 1 curve. For configurations of 1 curve in general
position one may ask if the lengths of all the arcs can be made equal (we
know that the answer is yes if the curve is orientation-reversing, and no
in general if the configuration is not in general position or contains more
than one curve). One may also ask if every minimal configuration of curves in
general position is contained in a configuration of shortest geodesics in
which each arc has the same length. These questions are equivalent to the
following:

\ques{Do all minimal 1-curve configurations
have close neighbors? Can every minimal configuration be extended to a
configuration with close neighbors? }

\begin{figure}[ht!]
\centerline{\includegraphics{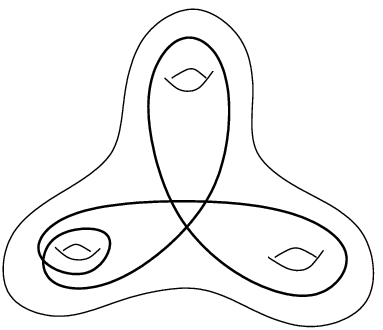}}
\nocolon\caption{}
\end{figure}

One may face strong restrictions when trying to extend a minimal
configuration to contain other curves. For example, if one wants to
extend the configuration in figure 12 to one containing a curve in the
homotopy class of figure 1, then the first curve must look as in figure 1b. 
Theorem 1.1 implies that some extension is always possible:

\begin{cor}Every minimal configuration of curves in $S$ can be extended to a
minimal configuration containing curves in any given homotopy classes in $S$.
\end{cor}

Denote by $l_{g}(a)$ the minimum length in the free homotopy class
of the curve $a$ when $S$ is given a Riemannian metric $g$. Denote by $%
a\bigcap b$ the minimum number of intersections between curves in the free
homotopy classes of $a$ and $b$, and by $a\bigcap a$ the minimal
number of self-intersections in the homotopy class of $a$.

\cor{If $l_{g}(a)\leq k\cdot l_{g}(b)$
for every Riemannian metric $g$ on $S$, then $a\bigcap c\leq k\cdot b\bigcap
c$ for every curve $c$ in $S$. In particular, $a\bigcap a\leq k\cdot
a\bigcap b\leq k^{2}\cdot b\bigcap b.$}

\proof Suppose that $a\bigcap c>k\cdot b\bigcap c$ for
some curve $c$. We may assume that $a$, $b$ and $c$ are in general position
and have minimal intersection and self-intersection. Apply the proof of
theorem 1.1 to the configuration formed by $a$ and $b$, but add to the
resulting measuring collection a copy of the curve $c$ with weight $w$. If
\[
w>\frac{k\cdot b\bigcap b+\left( k-1\right) \cdot a\bigcap b-a\bigcap a}{a\bigcap
c-k\cdot b\bigcap c} 
\]
then for the resulting metric $g$ we have
\[
{l_{g}(a)=a\bigcap a+a\bigcap b+w\cdot a\bigcap c>k\cdot }\left( {b\bigcap
a+b\bigcap b+w\cdot b\bigcap c}\right) ={k\cdot l_{g}(b)} 
\]
contrary to the hypothesis that $l_{g}(a)\leq k\cdot l_{g}(b)$. \endproof

Corollary 1.6 clearly holds if the curves $a$ , $b$ and $c$ are replaced by
any finite families of curves or arcs.

\section{Cutting and pasting.}

Let $\left\{ a_{i}\right\} $ be a configuration of curves with transverse
intersections in $S$. A {\sl cut and paste} on $\left\{a_{i}\right\} $ is 
done by cutting these curves at some of their intersection points and
glueing the resulting arcs in a different order to obtain a new collection of
curves $\left\{ b_{j}\right\}$. These curves have some `corners' that can
be rounded so the total number of intersections and the total length
of the original configuration are reduced.

\lem{If a collection of curves $\left\{ {a}_{i}\right\}$ in ${S}$ can be cut 
and pasted to obtain the collection $\left\{ {b}_{j}\right\}$, and
$\left\{ {a}_{i}\right\}$ can be homotoped to a collection
$\left\{ {a}_{i}^{\prime }\right\}$ without removing any
intersection points in the process, then $\left\{ {a}_{i}^{\prime }\right\}$
can be cut and pasted to obtain a collection homotopic to (the
nontrivial) $\left\{ {b}_{j}\right\}$.}

\proof We want to show that the instructions for cutting
and pasting $\left\{ {a}_{i}\right\} $ to get $\left\{ {b}_{j}\right\} $ can
be transmitted along the homotopy from $\left\{ {a}_{i}\right\} $ to $%
\left\{ {a}_{i}^{\prime }\right\} $ so that the final result is homotopic to 
$\left\{ {b}_{j}\right\} $. This is not obvious even though the intersection
points of $\left\{ {a}_{i}\right\} $ can be traced along the homotopy (they
don't disappear), because the result of doing the ``same'' cut and paste
before or after the homotopy may be different, as shown in figure 13.

\begin{figure}[ht!]
\centerline{\includegraphics{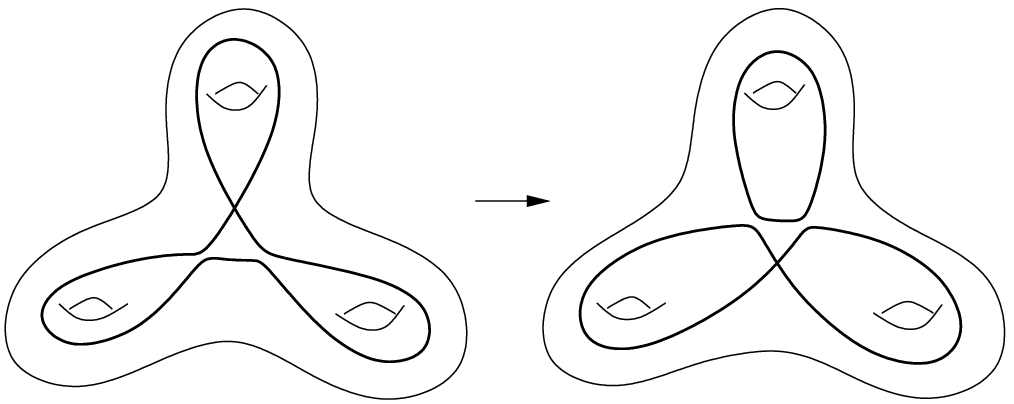}}
\nocolon\caption{}
\end{figure}

Any homotopy that doesn't remove intersection points can be done using 3
types of local moves in the configuration. The first two moves, adding a
small loop and creating a small bigon, do not change the homotopy class of
the resulting curves. Nevertheless, when doing these moves one can add cut
and paste instructions at the new intersections to avoid increasing the
number of intersections of the resulting curves (see figure 14a,b).

\begin{figure}[ht!]
\centerline{\includegraphics{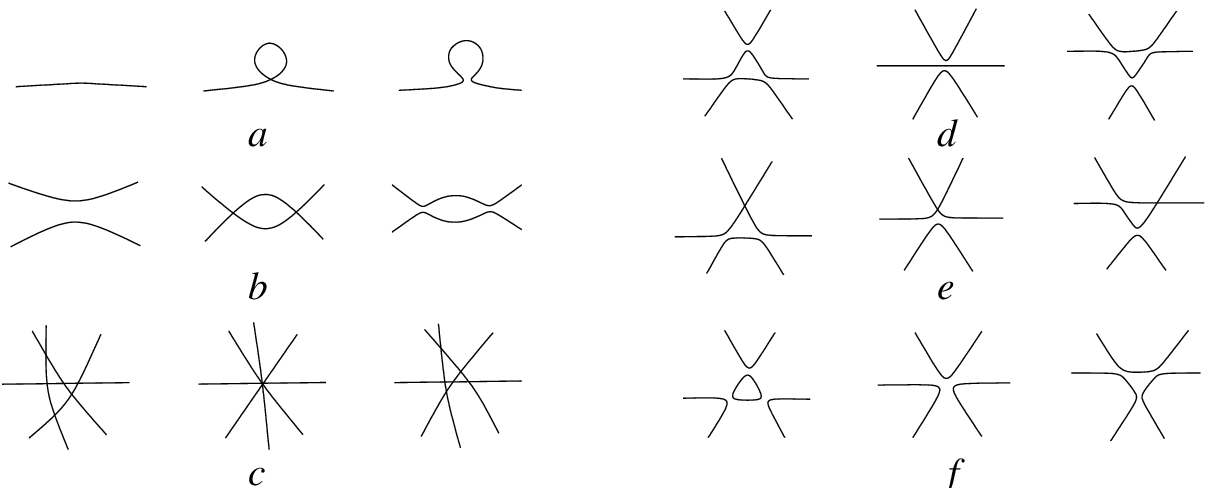}}
\nocolon\caption{}
\end{figure}

In the third move a local configuration of $n$ arcs that intersect each
other at different points collapses into one where all the arcs meet at a
single point, or viceversa: a configuration with a multiple intersection
opens up (see figure 14c). Observe that to transmit some cut and paste
instructions during these local moves one only needs that the endpoints of
the arcs that were connected by the original cut and paste instructions get
connected by the new instructions (any curve contained in the local
configuration is trivial). As the endpoints of the arcs that meet at a
single point can be connected at will by cutting and pasting at that point,
then all cut and paste instructions can be transmitted when a local
configuration collapses into a multiple intersection.

So the problem is to transmit the cut and paste instructions when a multiple
intersection opens up. In the case $n=3$ one can see how this can be done
directly (figure 14d-f shows some cases). Observe that the new instructions
may not be unique, but they can always be chosen to avoid creating new
curves and to avoid increasing the number of intersections of the resulting
curves.
In the case $n>3$, modify the homotopy so the multiple
intersection opens up one arc at a time. If an arc $a$ moves away from the
multiple intersection point and the cut and paste instructions don't change,
then the only connections that are affected are those involving the
endpoints of $a$, which are connected to the endpoints of at most 2 other
arcs of the local configuration. But we can change the cut and paste
instructions at the intersections of these 3 arcs as in the case $n=3$ to
get the right connections for the endpoints of $a$, and then change the
cutting and pasting instructions at the multiple intersection point as
needed to get the right connections between all the other endpoints. Now
repeat the argument until the multiple intersection opens up
completely.\endproof

In \cite{HS} Hass and Scott defined a `curve flow' that takes any
configuration of primitive curves in a surface to a configuration of
shortest geodesics by a homotopy that does not increase the number of
intersections at any moment. This result and the previous lemma imply the
following version of the theorem of Freedman, Hass and Scott:

\pro{Any finite family of primitive,
orientation-preserving curves in $S$ can be cut and pasted to obtain a
freely homotopic family of curves with minimal intersections and
self-intersections.}

\proof By \cite{HS} there is a homotopy that takes the
family $\left\{ {a}_{i}\right\} $ to some minimal intersection family $%
\left\{ {a}_{i}^{\prime }\right\} $ without increasing the number of
intersections, so running the homotopy backwards we get a homotopy that
takes $\left\{ {a}_{i}^{\prime }\right\} $ to $\left\{ {a}_{i}\right\} $
without removing any intersection points. Now lemma 2.1 shows how to transmit
the "don't cut anything" instructions in $\left\{ {a}_{i}^{\prime
}\right\} $ to cutting and pasting instructions in $\left\{ {a}_{i}\right\} $
without increasing the number of intersections of the resulting
curves.\endproof

Figure 15 shows a nonminimal configuration of 2 curves and a cut and paste
that transforms it into a minimal configuration.

\begin{figure}[ht!]
\centerline{\includegraphics{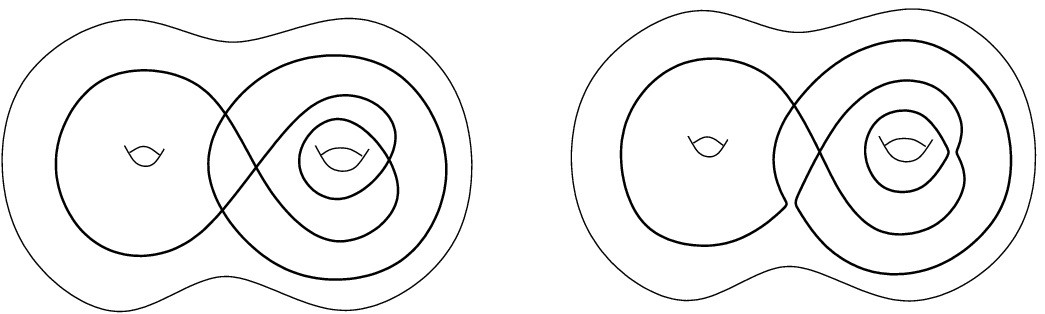}}
\nocolon\caption{}
\end{figure}

\cor{If a collection $\left\{a_{i}\right\} $ of curves with minimal intersection 
and self-intersection in $S$ can be cut and pasted to obtain the collection
$\left\{ b_{j}\right\}$, then $l_{g}(\left\{ {a}_{i}\right\} )>l_{g}(\left\{
b_{j}\right\} )$ for every Riemannian metric $g$ on $S$.}

\proof Observe that the hypothesis that $\left\{ {a}_{i}\right\} $ 
has minimal intersections is essential.
If $g$ is a Riemannian metric on $S$ and $\left\{ {a}_{i}^{\prime }\right\} $ is 
a collection of shortest geodesics (for the
metric $g$) homotopic to $\left\{ {a}_{i}\right\} $, then by \cite{HS} there
is a homotopy from $\left\{ {a}_{i}\right\} $ to $\left\{ {a}_{i}^{\prime
}\right\} $ that does not increase the number of intersections, so as $%
\left\{ {a}_{i}\right\} $ already had minimal intersections the number of
intersections must remain constant. So by lemma 2.1 the cutting and pasting
instructions to get $\left\{ b_{j}\right\} $ from $\left\{ {a}_{i}\right\} $
can be transmitted to get a homotopic collection $\left\{ b'_{j}\right\} $ from 
$\left\{ {a}_{i}^{\prime }\right\} $, so $l_{g}(\left\{ {a}%
_{i}\right\} )=l_{g}(\left\{ {a}_{i}^{\prime }\right\} )>l_{g}(\left\{
b'_{j}\right\} )=l_{g}(\left\{ b_{j}\right\} )$.\endproof

\exa{The converse to corollary 2.3 is not true.
Figure 16 shows 2 curves $a$ and $b$ on a surface such that $a$ cannot
be cut and pasted to obtain a curve homotopic to $b$ but one can
show that $l_{g}(a)>l_{g}(b)$ for every Riemannian metric $g$ on $S$ 
(so $a\bigcap c\geq b\bigcap c$ for every curve $c$).}

\begin{figure}[ht!]
\centerline{\includegraphics{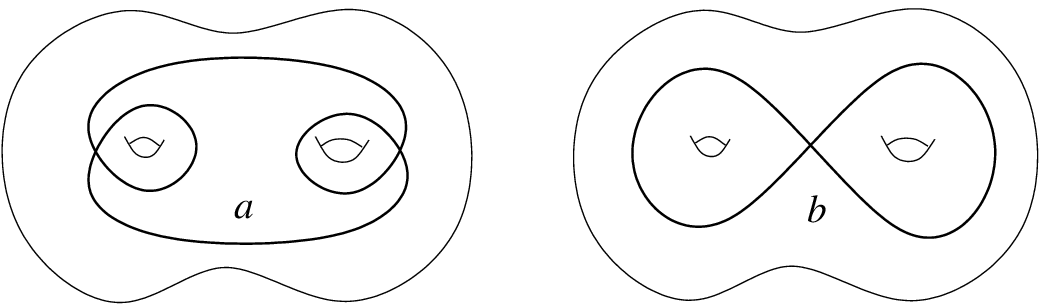}}
\nocolon\caption{}
\end{figure}

\que{If $a \bigcap c \leq b \bigcap c$ for
every curve $c$ in $S$, is it true that $l_g(a) \leq l_g(b)$ for every
Riemannian metric $g$ on $S$?}

\Addresses\recd 

\end{document}